\documentclass[reqno,12pt]{amsart}
\usepackage{amsfonts}
\usepackage{amssymb}
\usepackage{bbm}
\usepackage{times}

\def\B{\mathcal{B}}
\def\D{\mathbb{D}}

\def\C{\mathbb{C}}

\def\N{\mathbb N}

\def\B{\mathcal{B}}

\def\msk{\medskip}

\def\bege{\begin{equation}} \def\ende{\end{equation}}

\def\begr{\begin{eqnarray}} \def\endr{\end{eqnarray}}

\def\bege{\begin{equation}} \def\ende{\end{equation}}
\def\begr{\begin{eqnarray}} \def\endr{\end{eqnarray}}
\def\bnum{\begin{enumerate}} \def\enum{\end{enumerate}}

\begin{document}

\title[Differences of Weighted Composition Operators
]{ Differences of Weighted Composition Operators
from $H^\infty$ to the Bloch Space }
\author{Yecheng Shi,  Songxiao Li$^*$ and Xiangling Zhu}

\address{Yecheng Shi, Faculty of Information Technology,
Macau University of Science and Technology, Avenida Wai Long,
Taipa, Macau.
}\email{ 09ycshi@sina.cn}

\address{Songxiao Li, Institute of Fundamental and Frontier Science, University of Electronic Science and Technology of China, Avenue Xiyuan,
Chengdu, Sichuan, P.R.China.\newline
Faculty of Information Technology,
Macau University of Science and Technology, Avenida Wai Long,
Taipa, Macau.
}\email{jyulsx@163.com}

\address{Xiangling Zhu, Department of Mathematics, Jiaying University,  Meizhou 514015, China.  } \email{jyuzxl@163.com}

\begin{abstract}   In this paper, we give two new characterizations for the boundedness and compactness of the difference of two weighted composition
operators acting from $H^\infty$ to the Bloch space.
\thanks{*Corresponding author.}
\end{abstract}

\maketitle

\section{Introduction}

Let $\D$ be the open unit disk in the complex plane $\C$,   and $H(\D)$ be the class of functions analytic in $\D$.
We denote by $S(\D)$ the set of all analytic self-maps of $\D$.
For $a\in\D$, let $\sigma_{a}(z):=\frac{a-z}{1-\bar{a}z}$ be the automorphism of $\D$ exchanging 0 for $a$. For $z$, $w\in\D$, the pseudo-hyperbolic distance between $z$ and $w$ is given by
$$\rho(z,w)=|\sigma_{w}(z)|=\bigg|\frac{z-w}{1-\bar{w}z}\bigg|.$$
It is well known that $\rho(z,w)\leq1$ and
$$\frac{(1-|z|^2)(1-|w|^2)}{|1-\overline{w}z|^2} = 1-\rho(z,w)^2.$$
We write that $\rho(z)=\rho(\varphi(z),\psi(z))$ for the pseudo-hyperbolic distance between
$\varphi(z)$ and $\psi(z)$.

For a function $u\in H(\D)$ and $\varphi \in S(\D)$,  we define a weighted composition
operator $uC_{\varphi}$ on $H(\D)$ by
$$
(uC_\varphi f)(z) = u(z)f( \varphi(z)) ~~~~ ~\\~\\~\\~\mbox{~~~~for}~~f \in H(\D).
$$
When $u\equiv1$, we just have the composition operator $C_\varphi$ defined by $C_\varphi(f)=f\circ\varphi$.

Let $H^\infty = H^\infty(\D)$ be the space of all bounded analytic functions on $\D$.
Then $H^\infty$ is a Banach algebra with the supremum norm
$$\|f\|_{\infty} = \sup\{|f(z)|: z\in\D\}.$$

Recall that the Bloch space
$\mathcal{B}$, is a space which consists of all $f\in H(\mathbb{D})$ such that
$$
\|f\|_{\B}=\sup_{z \in \mathbb{D}}(1-|z|^2) |f'(z)| < \infty.
$$
It is well known that $\B$ is a Banach space under the norm $ |\|f\||_{\B}=|f(0)|+\|f\|_{\B}$.

We introduce some notations which will be used in the sequel. For $u\in H(\D)$ and $\varphi \in S(\D)$, let $\bigtriangleup$ be the set of all convergent sequence $\{z_k\}$ in $\D$.
Put
$$\Lambda_{u,\varphi} = \left\{\{z_k\} \in \triangle : |\varphi(z_k)| \to 1, (1-|z_k|^2)|u^{\prime}(z_k)| \nrightarrow 0\right\}$$
and
$$\Gamma^{\#}_{u,\varphi} = \left\{\{z_k\} \in \triangle : |\varphi(z_k)| \to 1, |u(z_k)\varphi^{\#}(z_k)| \nrightarrow 0\right\},$$
where $\varphi^{\#}(z)=\frac{1-|z|^2}{1-|\varphi(z)|^2}\varphi^{\prime}(z)$.
Then $uC_{\varphi}$ from $H^\infty$ to $\B$  is compact if and only if $\Lambda_{u,\varphi} = \Gamma^{\#}_{u,\varphi} = \emptyset$ (see \cite{O1}).\msk

Recently, many researcher  have studied the differences of composition operators, as well as the differences of weighted composition operators on some analytic function spaces. The main purpose for the study of the differences of composition operators is to
understand the topological structure of the set of composition
operators acting on the given function space. This line of research was first started in the setting of Hardy spaces, see \cite{Be,Sj}.
After that, such related problems have been studied on several classical function spaces such as $H^\infty$ \cite{HO3, MOZ}, the Bloch type spaces \cite{HO4, HO,li12, Ni,SL, SL2}, weighted-type spaces \cite{HLS,Ni}, the weighted Bergman spaces \cite{AW,Mo, S1,S2}, and the Hardy spaces \cite{G,S1}.

In \cite{wzz}, Wulan, Zheng and Zhu proved that $C_\varphi:\B\to\B$ is compact if and only if
$$\lim_{n\to\infty}\|\varphi^n\|_{\B}=0.$$  In \cite{SL}, we proved that $C_\varphi-C_\psi:\B\to\B$ is compact if and only if
$$\lim_{n\to\infty}\|\varphi^n-\psi^n\|_{\B}=0.$$  In \cite{HO1}, Hosokawa and Ohno studied the boundedness and compactness of the differences of two weighted composition operators acting from $H^\infty$ to
$\B$, under the following assumptions:
$$(A)~~U=\sup_{z\in\D}(1-|z|^2)|u'(z)|\rho(z) < \infty~~\\~\\~\mbox{and}~~\\~\\~V=\sup_{z\in\D}(1-|z|^2)|v'(z)|\rho(z) < \infty,$$
and
$$(B)~~~~(1-|z|^2)|u'(z)|\rho(z) \to 0~~\\~\\~\mbox{and}~~\\~\\~(1-|z|^2)|v'(z)|\rho(z) \to 0~~~\\~~\\~\mbox{if}~~\rho(z) \to 0.$$\msk

The main results in \cite{HO1} are stated as follows.\msk

\noindent{\bf Theorem A. \cite{HO1}} {\it Let u,v be in $H(\D)$ and $\varphi$, $\psi$ be in $S(\D)$. Assume that condition (A) holds. Then
$uC_{\varphi}-vC_{\psi}$ is bounded from $H^\infty$ to $\B$ if and only if the following conditions hold:\msk

(i) $u-v\in\B.$

(ii) $\sup_{z\in\D}|u(z)\varphi^{\#}(z)-v(z)\psi^{\#}(z)| < \infty.$

(iii) $\sup_{z\in\D}|v(z)\psi^{\#}(z)|\rho(z) < \infty.$\msk

 It is possible to replace $\varphi$ with $\psi$ and $u$ with $v$ in condition $(iii)$.
}\msk

\noindent{\bf Theorem B. \cite{HO1}} {\it Let u,v be in $H(\D)$ and $\varphi$, $\psi$ be in $S(\D)$. Assume that conditions (A) and (B) hold.
Suppose that $uC_{\varphi}-vC_{\psi}$ is bounded from $H^\infty$ to $\B$, but neither $uC_{\varphi}$ nor $vC_{\psi}$ is compact.
Then $uC_{\varphi}-vC_{\psi}$ is compact from $H^\infty$ to $\B$ if and only if the following conditions hold:

(i) $\Lambda_{u,\varphi} = \Lambda_{v,\psi}$.

(ii) $\lim_{n \to \infty}(1-|z_n|^2)|u'(z_n)-v'(z_n)| = 0$~~ for any $\{z_n\} \in \Lambda_{u,\varphi}$.

(iii) $\lim_{n \to \infty}(1-|z_n|^2)|u'(z_n))|\rho(z_n) = 0$~~ for any $\{z_n\} \in \Lambda_{u,\varphi}$.

(iv) $\Gamma^{\#}_{u,\varphi} = \Gamma^{\#}_{v,\psi}$.

(v) $\lim_{n \to \infty}|u(z_n)\varphi^{\#}(z_n)-v(z_n)\psi^{\#}(z_n)| = 0$~~ for any $\{z_n\} \in \Gamma^{\#}_{u,\varphi}$.

(vi) $\lim_{n \to \infty}|u(z_n)\varphi^{\#}(z_n)|\rho(z_n) = 0$~~ for any $\{z_n\} \in \Gamma^{\#}_{u,\varphi}$.\msk

It is possible to replace $u$ with $v$ in condition $(iii)$, and also
$u$ with $v$, and $\varphi$ with $\psi$ in condition $(vi)$.
}\msk

 Motivated by \cite{HO1,SL}, in this paper we will give some new characterizations for the boundedness and compactness of the differences of  weighted composition operators acting from $H^\infty$ to $\B$.
 More precisely, we show that $uC_\varphi-vC_\psi:H^\infty\to \B$ is bounded (respectively, compact) if and only if $$\sup_{n\in\N}\|u\varphi^n-v\psi^n\|_{\B}<\infty$$ (respectively,
 $\lim_{n\to\infty}\|u\varphi^n-v\psi^n\|_{\B}=0$), under some mild conditions.

For two quantities $A$ and $B$, the symbol $A\approx B$ means that $A\lesssim B\lesssim A$. We say that $A\lesssim B$ if there exists a constant $C$ such that $A\leq CB$.

\section{Main results and proofs}\vspace{0.2truecm}
In this section, we will give the main results and the proofs of this paper. We split this section into two subsection.

\subsection{ The boundedness of $uC_{\varphi}-vC_{\psi}: H^\infty\to \B$.}\vspace{0.2truecm}
In this subsection, we will give some new characterizations for the boundedness of the operator $uC_\varphi-vC_\psi:H^\infty\to\B$.
For this purpose we need the following lemmas.\msk

\noindent{\bf Lemma 2.1.} \cite{HO1} {\it For any $z,w\in\D$,
$$\bigg|1-\frac{(1-|z|^2)(1-|w^2|)}{(1-\overline{w}z)^2}\bigg| \leq 3\rho(z,w).$$}

\noindent{\bf Lemma 2.2.} \cite{HO} {\it For any $z,w\in\D$, we have $\rho(z,w)^2 \leq \flat(z,w) \leq 18\rho(z,w), $
where $$\flat(z,w)=\sup_{\|f\|_{\B\leq1}}|(1-|z|^2) f'(z)-(1-|w^2|) f'(w)|. $$}

The following lemma is basically from the proof of Theorem 3.1 of \cite{HO1}. \msk

\noindent{\bf Lemma 2.3.}  {\it Let $u,v$ be in $H(\D)$ and $\varphi$, $\psi$ be in $S(\D)$.   Then the following inequalities hold:
\begr
(i)&~&\sup_{z\in\D}|v(z)\psi^{\#}(z)|\rho(z)\nonumber\\
&~\leq& U+
\sup_{a\in\D}\|(u\varphi-v\psi)\sigma_{a}^2\|_{\B}
+\sup_{(a,b)\in\D^2}\|(u\varphi-v\psi)\sigma_{a}^2\sigma_{b}\|_{\B}.\nonumber
\endr
\begr
(ii)&~&\sup_{z\in\D}|u(z)\varphi^{\#}(z)-v(z)\psi^{\#}(z)|\nonumber\\
&~\leq& V+3U+\sup_{a\in\D}\|(u\varphi-v\psi)\sigma_{a}\|_{\B}+\nonumber\\
&&3\sup_{a\in\D}\|(u\varphi-v\psi)\sigma_{a}^2\|_{\B}+
3\sup_{(a,b)\in\D^2}\|(u\varphi-v\psi)\sigma_{a}^2\sigma_{b}\|_{\B}.\nonumber
\endr
}
 The following lemma is key for our first result in this paper. It has independent of interesting. \msk

\noindent{\bf Lemma 2.4.}  {\it Let u,v be in $H(\D)$ and $\varphi$, $\psi$ be in $S(\D)$. Let $n, m \in \N$. Then
$$\sup_{(a,b)\in\D^2}\|(uC_{\varphi}-vC_{\psi})\sigma_a^n\sigma_b^m\|_{\B} \lesssim
\sup_{k\in\N}\|u\varphi^k-v\psi^k\|_{\B}.$$}

{\it Proof.} Note that the Maclaurin expansion of the function $\sigma^n_a\sigma^m_b$ is given by
\begr
&&(\sigma_a^n\sigma_b^m)(z)\nonumber \\
&=&\left(a-(1-|a|^2)\sum_{k=0}^{\infty}\overline{a}^kz^{k+1}\right)^n
\left(b-(1-|b|^2)\sum_{k=0}^{\infty}\overline{b}^kz^{k+1}\right)^m\nonumber \\
&=&\left(a^n+\sum_{l=1}^{n}C_{n}^{l}a^{n-l}(|a|^2-1)^l\sum_{k=l}^{\infty}C_{k-1}^{k-l}\overline{a}^{k-l}z^k\right)\nonumber \\
&~~&\cdot\left(b^m+\sum_{r=1}^{m}C_{m}^{r}b^{m-r}(|b|^2-1)^r\sum_{k=r}^{\infty}C_{k-1}^{k-r}\overline{b}^{k-r}z^k\right)
\nonumber
\endr
\begr
&=&\sum_{l=1}^{n}\sum_{r=1}^{m}\bigg\{C_{n}^{l}C_{m}^{r}a^{n-l}b^{m-r}(|a|^2-1)^l(|b|^2-1)^r\nonumber \\
&~~&\cdot(\sum_{k=l}^{\infty}C_{k-1}^{k-l}\overline{a}^{k-l}z^k)
(\sum_{k=r}^{\infty}C_{k-1}^{k-r}\overline{b}^{k-r}z^k)\bigg\}\nonumber \\
&~&+a^n\left(\sum_{r=1}^{m}C_{m}^{r}b^{m-r}(|b|^2-1)^r\sum_{k=r}^{\infty}C_{k-1}^{k-r}\overline{b}^{k-r}z^k\right)\nonumber \\
&~&+b^m\left(\sum_{l=1}^{n}C_{n}^{l}a^{n-l}(|a|^2-1)^l\sum_{k=l}^{\infty}C_{k-1}^{k-l}\overline{a}^{k-l}z^k\right)+a^nb^m\nonumber\\
&=&\sum_{l=1}^{n}\sum_{r=1}^{m}\bigg\{C_{n}^{l}C_{m}^{r}a^{n-l}b^{m-r}(|a|^2-1)^l(|b|^2-1)^r\nonumber \\
&~~&\cdot\sum_{k=l+r}^{\infty}\left(\sum_{s=l}^{k-r}C_{s-1}^{s-l}C_{k-s-1}^{k-r-s}\overline{a}^{s-l}\overline{b}^{k-r-s}\right)z^k\bigg\}
\nonumber\\
&~&+a^n\left(\sum_{r=1}^{m}C_{m}^{r}b^{m-r}(|b|^2-1)^r\sum_{k=r}^{\infty}C_{k-1}^{k-r}\overline{b}^{k-r}z^k\right)\nonumber \\
&~&+b^m\left(\sum_{l=1}^{n}C_{n}^{l}a^{n-l}(|a|^2-1)^l\sum_{k=l}^{\infty}C_{k-1}^{k-l}\overline{a}^{k-l}z^k\right)+a^nb^m\nonumber \\
&:=&f_1+f_2+f_3+f_4,\nonumber
\endr
where we used the fact that
\begr
(\sum_{k=0}^\infty\overline{a}^kz^{k+1})^l=\sum_{k=l}^\infty C^{k-l}_{k-1}\overline{a}^{k-l}z^k
~~\mbox{~~and~~} ~~
(\sum_{k=0}^\infty\overline{b}^kz^{k+1})^r=\sum_{k=l}^\infty C^{k-r}_{k-1}\overline{b}^{k-r}z^k.\nonumber
\endr
It follows from the triangle inequality  that
\begr
&&\|(uC_{\varphi}-vC_{\psi})\sigma_a^n\sigma_b^m\|_{\B}\nonumber\\
&\leq&\|(uC_{\varphi}-vC_{\psi})f_1\|_{\B}+\|(uC_{\varphi}-vC_{\psi})f_2\|_{\B}\nonumber\\
& & +\|(uC_{\varphi}-vC_{\psi})f_3\|_{\B}+\|(uC_{\varphi}-vC_{\psi})f_4\|_{\B}.\nonumber
\endr
For $l\geq1$, $r\geq1,$ after a calculation, we have \begr
\sum_{k=l+r}^{\infty}\sum_{s=l}^{k-r}C_{s-1}^{s-l}C_{k-s-1}^{k-s-r}|a|^{s-l}|b|^{k-s-r}
&=&\sum_{n=0}^{\infty}\sum_{m=0}^{\infty}C_{n+l-1}^{n}C_{m+r-1}^{m}|a|^{n}|b|^{m}\nonumber\\
&\lesssim& \sum_{n=0}^{\infty}n^{l-1}|a|^{n}\sum_{m=0}^{\infty}m^{r-1}|b|^{m}\nonumber\\
&\lesssim&(1-|a|^2)^{-l}(1-|b|^2)^{-r}.\nonumber
\endr
Here, we used $C_{n+l-1}^{n}\approx n^{l-1}, ~~~\mbox{as}~~n\to\infty.$
Therefore,
\begr
&&\|(uC_{\varphi}-vC_{\psi})f_1\|_{\B}\nonumber\\
&\leq&\sum_{l=1}^{n}\sum_{r=1}^{m}C_{n}^{l}C_{m}^{r}|a|^{n-l}|b|^{m-r}(1-|a|^2)^l(1-|b|^2)^r\nonumber\\
&&\sum_{k=l+r}^{\infty}\bigg(\sum_{s=l}^{k-r}C_{s-1}^{s-l}C_{k-s-1}^{k-s-r}|a|^{s-l}|b|^{k-s-r}\bigg)\|u\varphi^k-v\psi^k\|_{\B}\nonumber\\
&\lesssim&\sum_{l=1}^{n}\sum_{r=1}^{m}C_{n}^{l}C_{m}^{r}|a|^{n-l}|b|^{m-r}(1-|a|^2)^l(1-|b|^2)^r\nonumber\\
&&\sum_{k=l+r}^{\infty}\bigg(\sum_{s=l}^{k-r}C_{s-1}^{s-l}C_{k-s-1}^{k-s-r}|a|^{s-l}|b|^{k-s-r}\bigg)
\sup_{n\in\N}\|u\varphi^n-v\psi^n\|_{\B}\nonumber\\
&\lesssim&\sum_{l=1}^{n}\sum_{r=1}^{m}C_{n}^{l}C_{m}^{r}|a|^{n-l}|b|^{m-r}\sup_{n\in\N}\|u\varphi^n-v\psi^n\|_{\B}\nonumber\\
&\lesssim&\sup_{n\in\N}\|u\varphi^n-v\psi^n\|_{\B},\nonumber
\endr
where we used the fact that
$$\sum_{l=1}^{n}\sum_{r=1}^{m}C_{n}^{l}C_{m}^{r}|a|^{n-l}|b|^{m-r}\leq(1+|a|)^n(1+|b|)^m\leq2^{n+m}.$$
Since
$$\sum_{k=r}^{\infty}C_{k-1}^{k-r}|b|^{k-r}\lesssim\sum_{k=r}^{\infty}k^{r-1}|b|^{k-r}\lesssim\sum_{k=0}^{\infty}k^{r-1}|b|^{k}\lesssim(1-|b|^2)^{-r},$$
we have
\begr
&&|a|^n\left(\sum_{r=1}^{m}C_{m}^{r}|b|^{m-r}(1-|b|^2)^r\sum_{k=r}^{\infty}C_{k-1}^{k-r}|b|^{k-r}\right)
\lesssim \sum_{r=1}^{m}C_{m}^{r}|b|^{m-r}
\lesssim 1.\nonumber
\endr
Using the above estimates and the triangle inequality, we get
\begr
&&\|(uC_{\varphi}-vC_{\psi})f_2\|_{\B}\nonumber\\
&\leq&|a|^n\sum_{r=1}^{m}C_{m}^{r}|b|^{m-r}(1-|b|^2)^r\sum_{k=r}^{\infty}C_{k-1}^{k-r}|b|^{k-r}\|u\varphi^k-v\psi^k\|_{\B}
\nonumber\\
&\lesssim&\sup_{n\in\N}\|u\varphi^n-v\psi^n\|_{\B}.\nonumber
\endr
Similarly, we obtain
\begr
\|(uC_{\varphi}-vC_{\psi})f_3\|_{\B}
\lesssim\sup_{n\in\N}\|u\varphi^n-v\psi^n\|_{\B}.\nonumber
\endr
In addition,
\begr
\|(uC_{\varphi}-vC_{\psi})f_4\|_{\B}\leq \|(uC_{\varphi}-vC_{\psi})1\|_{\B}
\leq\sup_{n\in\N}\|u\varphi^n-v\psi^n\|_{\B}.\nonumber
\endr
The proof is complete.\msk

Now we are in a position to state our first main result in this paper.\msk

\noindent{\bf Theorem 2.1.} {\it Let $u,v$ be in $H(\D)$ and $\varphi$, $\psi$ be in $S(\D)$. Assume that condition (A) holds. Then
the following statements are equivalent.

(i)~$uC_{\varphi}-vC_{\psi} : H^\infty \to \B$ is bounded.

(ii) For $(n,m) \in \{(0,0), (0,1), (0,2), (1,2)\}$, we have
$$\sup_{(a,b)\in\D^2}
\|(uC_{\varphi}-vC_{\psi})\sigma_a^n\sigma_b^m\|_{\B}<\infty.$$

(iii) $$\sup_{n\in\N} \|u\varphi^n-v\psi^n\|_{\B}<\infty.$$ }

\noindent{\it Proof.} $(i) \Rightarrow (iii)$. First, we assume that  $uC_{\varphi}-vC_{\psi}:H^\infty\to\B$ is bounded. Consider the function $f_n = z^n$, we have $\|f_n\|_{\infty} = 1$.
Thus, for any $k\in\N$,
$$\infty>\|uC_{\varphi}-vC_{\psi}\|_{H^\infty\to\B} \geq \sup_{k\in\N}\|(uC_{\varphi}-vC_{\psi})f_k\|_{\B}=\sup_{k\in\N}\|u\varphi^k-v\psi^k\|_{\B},$$
which implies the desire result.

$(iii)\Rightarrow (ii)$. By Lemma 2.4, we get it immediately.

$(ii)\Rightarrow (i)$. Let $f\in H^\infty$ such that $\|f\|_{\infty} \leq 1$. Then 
\begr
&&(1-|z|^2)\big((uC_{\varphi}-vC_{\psi})f\big)^{\prime}(z)\nonumber\\
&=&(1-|z|^2)\Big(u^{\prime}(z)f(\varphi(z))-v^{\prime}(z)f(\psi(z))\Big)\nonumber\\
&~&+(1-|z|^2)\Big(u(z)\varphi'(z)f'(\varphi(z))-v(z)\psi'(z)f'(\psi(z))\Big)\nonumber\\
&=&(1-|z|^2)(u'(z)-v'(z))f(\varphi(z))\nonumber\\
&~&+(1-|z|^2)v^{\prime}(z)\big(f(\varphi(z))-f(\psi(z))\big)\nonumber\\
&~&+\big(u(z)\varphi^{\#}(z)-v(z)\psi^{\#}(z)\big)(1-|\varphi(z)|^2)f'(\varphi(z))\nonumber\\
&~&+v(z)\psi^{\#}(z)\Big((1-|\varphi(z)|^2)f'(\varphi(z))-(1-|\psi(z)|^2)f'(\psi(z))\Big).\nonumber
\endr
It is known that for $z, w \in \D$ and $f \in H^\infty$,
$$\sup_{\|f\|_{\infty}\leq1}|f(z)-f(w)| = \frac{2-2\sqrt{1-\rho(z,w)^2}}{\rho(z,w)} \leq 2\rho(z,w).$$
Thus, combining the above facts, Lemma 2.2 and condition (A), we obtain
\begr
&&\|uC_{\varphi}-vC_{\psi}\|_{H^\infty\to\B}\nonumber\\
&\leq&  |u(0)|+|v(0)|+\|u-v\|_{\B}+2V+ \sup_{z\in\D}|u(z)\varphi^{\#}(z)-v(z)\psi^{\#}(z)|+\nonumber\\
&~~& 18\sup_{z\in\D}|v(z)\psi^{\#}(z)|\rho(z). \nonumber
\endr
  Therefore, using Lemmas 2.3 and 2.4, we  obtain
\begr
&&\|uC_{\varphi}-vC_{\psi}\|_{H^\infty\to\B}\nonumber\\
&\leq&  |u(0)|+|v(0)|+21U+3V+\|u-v\|_{\B}+\sup_{a\in\D}\|(u\varphi-v\psi)\sigma_{a}\|_{\B}\nonumber\\
& &+21\sup_{a\in\D}\|(u\varphi-v\psi)\sigma_{a}^2\|_{\B}
+21\sup_{(a,b)\in\D^2}\|(u\varphi-v\psi)\sigma_{a}^2\sigma_{b}\|_{\B}\nonumber\\
& <&\infty.\nonumber
\endr
The proof is complete.

\subsection{The compactness of $uC_{\varphi}-vC_{\psi}:H^\infty\to\B$}\vspace{0.2truecm}
In this subsection, we will give two new characterizations for the compactness of the differences of two weighted composition operators acting
from $H^\infty$ to $\B$.\msk

\noindent{\bf Lemma 2.5.}  {\it Let $u,v$ be in $H(\D)$ and $\varphi$, $\psi$ be in $S(\D)$. Let $n, m \in \N$.
Then
$$\lim_{|a| \to 1}\sup_{b\in\D}\|(uC_{\varphi}-vC_{\psi})(\sigma_a^{n+1}-a\sigma_a^n)\sigma_b^m\|_{\B} \lesssim
\limsup_{k \to \infty}\|u\varphi^k-v\psi^k\|_{\B}.$$ }

\noindent{\it Proof.}  Let $g=(\sigma_a^{n+1}-a\sigma_a^n)\sigma_b^m$.
As in Lemma 2.4,   we have
$$(\sigma_a^n\sigma_b^m)(z)=f_1(z)+f_2(z)+f_3(z)+f_4(z),$$
where $f_i(i=1,2,3,4) $ are defined in the proof of Lemma 2.4. Thus,
$$g(z)=(f_1+f_2+f_3+f_4)\cdot(\sigma_a(z)-a):=g_1+g_2+g_3+g_4,$$
where $g_i(z)=f_i\cdot(\sigma_a(z)-a), ~~~~i=1, 2, 3, 4.$  It is easy to see that
$$\sigma_a(z)-a = -(1-|a|^2)\sum_{k=0}^{\infty}\overline{a}^kz^{k+1}.$$
Denote that $c_k(l,r) = \sum_{s=l}^{k-r}C_{s-1}^{s-l}C_{k-s-1}^{k-s-r}
\overline{a}^{s-l}\overline{b}^{k-s-r}$. We have

  \begr
g_1(z)&=&f_1(z)\left(\sigma_a(z)-a\right)\nonumber \\
&=&(|a|^2-1)\sum_{l=1}^{n}\sum_{r=1}^{m}\bigg\{C_{n}^{l}C_{m}^{r}a^{n-l}b^{m-r}(|a|^2-1)^l(|b|^2-1)^r\nonumber \\
&~~&\cdot\sum_{t=l+r+1}^{\infty}
\left(\sum_{i=0}^{t-l-r-1}c_{t-1-i}(l,r)\overline{a}^i\right)z^t\bigg\}.\nonumber
\endr
Here, we used
$$(\sum_{k=l+r}^\infty c_k(l,r)z^k)(\sum_{k=0}^\infty\overline{a}^kz^{k+1})=\sum_{t=l+r+1}(\sum_{i=0}^{t-l-r-1}c_{t-1-i(l,r)}\overline{a}^{i})z^t.$$
Using the triangle inequality and by elementary  calculations, we get
\begr
&&\|(uC_{\varphi}-vC_{\psi})\left(f_1\cdot\left(\sigma_a-a\right)\right)\|_{\B}\nonumber\\
&\leq&(1-|a|^2)\sum_{l=1}^{n}\sum_{r=1}^{m}\bigg\{C_{n}^{l}C_{m}^{r}|a|^{n-l}|b|^{m-r}(1-|a|^2)^l(1-|b|^2)^r\nonumber \\
&~~&\sum_{t=l+r+1}^{\infty}\left|\sum_{i=0}^{t-l-r-1}c_{t-1-i}(l,r)\overline{a}^i\right|\|u\varphi^t-v\psi^t\|_{\B}\bigg\}\nonumber\\
&\leq&(1-|a|^2)\sum_{l=1}^{n}\sum_{r=1}^{m}\bigg\{C_{n}^{l}C_{m}^{r}|a|^{n-l}|b|^{m-r}(1-|a|^2)^l(1-|b|^2)^r\nonumber \\
&&\sum_{t=l+r+1}^{\infty}F(t,l,r,a,b)\|u\varphi^t-v\psi^t\|_{\B}\bigg\}.
\nonumber
\endr
Here and henceforth,
$$F(t,l,r,a,b)=\sum_{i=0}^{t-l-r-1}\sum_{s=l}^{t-1-i-r}C_{s-1}^{s-l}C_{t-1-i-s-1}^{t-1-i-s-r}|a|^{s+i-l}|b|^{t-1-i-s-r}
.$$
For each $N > n+m$, we split the above quantity into two terms
\begr
&& I_1:=(1-|a|^2)\sum_{l=1}^{n}\sum_{r=1}^{m}\bigg\{C_{n}^{l}C_{m}^{r}|a|^{n-l}|b|^{m-r}(1-|a|^2)^l(1-|b|^2)^r\nonumber \\
&&\sum_{t=l+r+1}^{N}F(t,l,r,a,b)\|u\varphi^t-v\psi^t\|_{\B}\bigg\}
\nonumber
\endr
and
\begr
&& I_2:=(1-|a|^2)\sum_{l=1}^{n}\sum_{r=1}^{m}\bigg\{C_{n}^{l}C_{m}^{r}|a|^{n-l}|b|^{m-r}(1-|a|^2)^l(1-|b|^2)^r\nonumber \\
&&\sum_{t=N+1}^{\infty}F(t,l,r,a,b)\|u\varphi^t-v\psi^t\|_{\B}\bigg\}.
\nonumber
\endr
For the first term, letting $|a| \to 1$, we have  $\lim_{|a| \to 1}\sup_{b\in\D} I_1= 0.$ For the second term, we have
 \begr
I_2&\leq&(1-|a|^2)\sum_{l=1}^{n}\sum_{r=1}^{m}\bigg\{C_{n}^{l}C_{m}^{r}|a|^{n-l}|b|^{m-r}(1-|a|^2)^l(1-|b|^2)^r\nonumber \\
&&\sum_{t=N+1}^{\infty}F(t,l,r,a,b)\bigg\}\sup_{k\geq N+1}\|u\varphi^k-v\psi^k\|_{\B}.
\nonumber
\endr
 For fixed $l, r$, since $N>n+m\geq l+r$, we get
  \begr
I_2&\leq&(1-|a|^2)\sum_{l=1}^{n}\sum_{r=1}^{m}\bigg\{C_{n}^{l}C_{m}^{r}|a|^{n-l}|b|^{m-r}(1-|a|^2)^l(1-|b|^2)^r\nonumber \\
&&\sum_{t=l+r+1}^{\infty}F(t,l,r,a,b)\bigg\}\sup_{k\geq N+1}\|u\varphi^k-v\psi^k\|_{\B}.
\nonumber
\endr
Since
\begr
&&\sum_{k=l+r}^{\infty}\left(\sum_{s=l}^{k-r}C_{s-1}^{s-l}C_{k-s-1}^{k-s-r}|a|^{s-l}|b|^{k-s-r}\right)\sum_{n=0}^{\infty}|a|^n\nonumber\\
&=&\sum_{k=l+r}^{\infty}\sum_{n=0}^{\infty}\sum_{s=l}^{k-r}C_{s-1}^{s-l}C_{k-s-1}^{k-s-r}|a|^{s+n-l}|b|^{k-s-r}\nonumber\\
&=& \sum_{t=l+r+1}^{\infty}\sum_{i=0}^{t-l-r-1}
\sum_{s=l}^{t-1-i-r}C_{s-1}^{s-l}C_{t-1-i-s-1}^{t-1-i-s-r}|a|^{s+i-l}|b|^{t-1-i-s-r}\nonumber \\
&=&\sum_{t=l+r+1}^{\infty} F(t,l,r,a,b)
\nonumber
\endr
and
\begr
&&\sum_{k=l+r}^{\infty}\left(\sum_{s=l}^{k-r}C_{s-1}^{s-l}C_{k-s-1}^{k-r-s}|a|^{s-l}|b|^{k-r-s}\right)\sum_{n=0}^{\infty}|a|^n\nonumber\\
&\lesssim& (1-|a|^2)^{-l}(1-|b|^2)^{-r}(1-|a|)^{-1}, \nonumber
\endr
we have
\begr
I_2 \lesssim\sum_{l=1}^{n}\sum_{r=1}^{m}C_{n}^{l}C_{m}^{r}|a|^{n-l}|b|^{m-r}
\sup_{k\geq N+1}\|u\varphi^k-v\psi^k\|_{\B}.\nonumber
\endr
Noting that
\begr
\sum_{l=1}^{n}\sum_{r=1}^{m}C_{n}^{l}C_{m}^{r}|a|^{n-l}|b|^{m-r}&=&\left(\sum_{l=1}^{n}C_{n}^{l}|a|^{n-l}\right)
\left(\sum_{r=1}^{m}C_{m}^{r}|b|^{m-r}\right)\nonumber\\
&\leq&(1+|a|)^n(1+|b|)^m\nonumber\\
&\leq&2^{n+m},\nonumber
\endr
we obtain
\begr
\lim_{|a| \to 1}\sup_{b\in\D} I_2\lesssim \sup_{k\geq N+1}\|u\varphi^k-v\psi^k\|_{\B}.
\nonumber
\endr
Therefore,
$$\lim_{|a| \to 1}\sup_{b\in\D}\|(uC_{\varphi}-vC_{\psi})g_1\|_{\B} \lesssim
\limsup_{k \to \infty}\|u\varphi^k-v\psi^k\|_{\B}.$$

Since
\begr
&&g_2\nonumber\\
&=&(|a|^2-1)a^n\Big(\sum_{r=1}^{m}C_{m}^{r}b^{m-r}(|b|^2-1)^r\sum_{k=r}^{\infty}C_{k-1}^{k-r}\overline{b}^{k-r}z^k\Big)
\Big(\sum_{i=0}^{\infty}\overline{a}^iz^{i+1}\Big)
\nonumber\\
&=&(|a|^2-1)a^n\sum_{r=1}^{m}C_{m}^{r}b^{m-r}(|b|^2-1)^r\Big(\sum_{k=r}^{\infty}C_{k-1}^{k-r}\overline{b}^{k-r}z^k
\sum_{i=0}^{\infty}\overline{a}^iz^{i+1}\Big)\nonumber\\
&=&(|a|^2-1)a^n\sum_{r=1}^{m}C_{m}^{r}b^{m-r}(|b|^2-1)^r\sum_{k=r+1}^{\infty}\sum_{i=1}^{k-r}C_{k-i-1}^{k-i-r}\overline{b}^{k-i-r}\overline{a}^{i-1}z^k,
\nonumber
\endr
by the triangle inequality, we obtain
\begr
&&
\|(uC_{\varphi}-vC_{\psi})g_2\|_{\B}\nonumber\\
&\leq&(1-|a|^2)|a|^n\sum_{r=1}^{m}C_{m}^{r}|b|^{m-r}(1-|b|^2)^r\nonumber\\
&~&\sum_{k=r+1}^{\infty}\sum_{i=1}^{k-r}C_{k-i-1}^{k-i-r}|b|^{k-i-r}|a|^{i-1}\|u\varphi^k-v\psi^k\|_{\B}.
\nonumber
\endr
For any fixed integer $N>r$, we have
\begr
&&\|(uC_{\varphi}-vC_{\psi})g_2\|_{\B}\nonumber\\
&\leq&(1-|a|^2)|a|^n\sum_{r=1}^{m}C_{m}^{r}|b|^{m-r}(1-|b|^2)^r\nonumber\\
&~&\sum_{k=r+1}^{N}\sum_{i=1}^{k-r}C_{k-i-1}^{k-i-r}|b|^{k-i-r}|a|^{i-1}\sup_{n\in\N}\|u\varphi^n-v\psi^n\|_{\B}\nonumber\\
&~&+(1-|a|^2)|a|^n\sum_{r=1}^{m}C_{m}^{r}|b|^{m-r}(1-|b|^2)^r\nonumber\\
&~&\sum_{k=N+1}^{\infty}\sum_{i=1}^{k-r}C_{k-i-1}^{k-i-r}|b|^{k-i-r}|a|^{i-1}\sup_{n\geq N+1}\|u\varphi^n-v\psi^n\|_{\B}
\nonumber\\
&:=&J_1+J_2.\nonumber
\endr
After a calculation,
\begr
\sum_{k=N+1}^{\infty}\sum_{i=1}^{k-r}C_{k-i-1}^{k-i-r}|b|^{k-i-r}|a|^{i-1}&\leq&
\sum_{k=r+1}^{\infty}\sum_{i=1}^{k-r}C_{k-i-1}^{k-i-r}|b|^{k-i-r}|a|^{i-1}\nonumber\\
&=&\sum_{k=r}^{\infty}C_{k-1}^{k-r}|b|^{k-r}\sum_{i=1}^{\infty}|a|^{i-1}\nonumber\\
&\lesssim& (1-|b|^2)^{-r}(1-|a|^2)^{-1}. \nonumber
\endr
Letting $|a|\to1$, we obtain $\lim_{|a| \to 1}\sup_{b\in\D} J_1= 0$ and
$$
\lim_{|a| \to 1}\sup_{b\in\D} J_2 \lesssim \sup_{k\geq N+1}\|u\varphi^k-v\psi^k\|_{\B}.$$
Hence
$$
\lim_{|a| \to 1}\sup_{b\in\D}\|(uC_{\varphi}-vC_{\psi})g_2\|_{\B} \nonumber\\
\lesssim \limsup_{k\to\infty}\|u\varphi^k-v\psi^k\|_{\B}.$$
Similarly, we can prove that
$$\lim_{|a| \to 1}\sup_{b\in\D}\|(uC_{\varphi}-vC_{\psi})g_3\|_{\B} \lesssim
\limsup_{k \to \infty}\|u\varphi^k-v\psi^k\|_{\B}$$
and
$$\lim_{|a| \to 1}\sup_{b\in\D}\|(uC_{\varphi}-vC_{\psi})g_4\|_{\B} \lesssim
\limsup_{k \to \infty}\|u\varphi^k-v\psi^k\|_{\B}.$$
Therefore,
$$\lim_{|a| \to 1}\sup_{b\in\D}\|(uC_{\varphi}-vC_{\psi})g\|_{\B} \lesssim
\limsup_{k \to \infty}\|u\varphi^k-v\psi^k\|_{\B}.$$
The proof is complete. \msk

\noindent{\bf Theorem 3.2.} {\it Let $u,v$ be in $H(\D)$ and $\varphi$, $\psi$ be in $S(\D)$. Assume that the conditions (A) and (B) hold. Suppose $uC_{\varphi}-vC_{\psi}$ is bounded from $H^\infty$ to $\B$. Then
the following statements are equivalent.

(i)~$uC_{\varphi}-vC_{\psi} : H^\infty \to \B$ is compact.

(ii) For $(n,m) \in \{(0,0), (0,1), (0,2), (1,2)\}$, we have
$$\lim_{|a| \to 1}\sup_{b \in \D}\|(uC_{\varphi}-vC_{\psi})\left((\sigma_a^{n+1}-a\sigma_a^{n})\sigma_b^m\right)\|_{\B} =0.$$

(iii) $$\lim_{k \to \infty} \|u\varphi^k-v\psi^k\|_{\B} =0.$$ }

\noindent{\it Proof.} $(i) \Rightarrow (iii)$. We assume that $uC_{\varphi}-vC_{\psi}$ from $H^\infty$ to $\B$ is compact.
Consider the function $f_k(z) = z^k$. We have
$\|f_{k}\|_{\infty}=1$ and $\{f_k\}$ converges to $0$ uniformly on every compact subset of $\D$. Thus, $\lim_{k\rightarrow\infty}\|(uC_{\varphi}-vC_{\psi})f_k\|_{\B}=0$, i.e.,
$$\lim_{k\rightarrow\infty}\|u\varphi^k-v\psi^k\|_{\B} = 0.$$

$(iii) \Rightarrow (ii). $ It is obvious from Lemma 2.5.

$(ii) \Rightarrow (i)$. We consider it in two cases.

{\it Case 1.}  One of the weighted composition operators $uC_{\varphi}$ and $vC_\psi$ is compact. We may assume that
 $vC_\psi$ is compact and only have to prove $(ii)$ implies that $uC_{\varphi}$ is compact from $H^\infty$ to $\B$.
For any $w\in\D$, define
$$f_{\varphi(w)}(z) = \sigma_{\varphi(w)}^2(z)-\varphi(w)\sigma_{\varphi(w)}(z), ~~~z\in \D.$$
We have
$$\|uC_{\varphi}f_{\varphi(w)}\|_{\B} \geq (1-|w|^2)|u(w)\frac{1}{1-|\varphi(w)|^2}\varphi'(w)\varphi(w)|.$$
Thus
\begr
\lim_{|a| \to 1}\|uC_{\varphi}f_{a}\|_{\B} \geq \lim_{|\varphi(w)| \to 1}\|uC_{\varphi}f_{\varphi(w)}\|_{\B}
 \geq \lim_{|\varphi(w)| \to 1}|u(w)\varphi^{\#}(w)|.\nonumber
 \endr
Also, set
$$g_{\varphi(w)}(z) = \sigma_{\varphi(w)}(z)-\varphi(w), ~~~z\in \D.  $$
We obtain
$$\|uC_{\varphi}g_{\varphi(w)}\|_{\B} \geq |(1-|w|^2)u'(w)\varphi(w)+u(w)\varphi^{\#}(w)|.$$
Therefore
\begr
\lim_{|a| \to 1}\|uC_{\varphi}f_{a}\|_{\B}+\lim_{|a| \to 1}\|uC_{\varphi}g_{a}\|_{\B}
\geq\lim_{|\varphi(w)| \to 1}(1-|w|^2)|u'(w)|.\nonumber
 \endr
 Note that
 \begr \lim_{|a| \to 1}\sup_{b \in \D}
\|(uC_{\varphi}-vC_{\psi})\left((\sigma_a-a)\sigma_b\right)\|_{\B} \geq
 \lim_{|a| \to 1}\|(uC_{\varphi}-vC_{\psi})(\sigma_a^2-a\sigma_a)\|_{\B}. \nonumber
\endr
By $(ii)$, we have
$$ \lim_{|\varphi(w)| \to 1}|u(w)\varphi^{\#}(w)|=\lim_{|\varphi(w)| \to 1}(1-|w|^2)|u'(w)|=0.$$
Hence,
by Theorem 3 of \cite{O1} we get that $uC_\varphi:H^\infty\to \B$ is compact, as desired.

 {\it Case 2.}  If neither $uC_{\varphi}$ nor $vC_\psi$ is compact.
 We need only to prove that condition $(ii)$ implies conditions $(i)$-$(vi)$ of Theorem B.  We assume that condition $(ii)$  holds.

Suppose that $\Gamma^{\#}_{u,\varphi} \neq \emptyset.$ For a sequence $\{z_n\}\in \Gamma^{\#}_{u,\varphi},$
we consider the function
$$f_{1,z_n}(z) = (\sigma^2_{\varphi(z_n)}(z)-\varphi(z_n)\sigma_{\varphi(z_n)}(z))\sigma_{\psi(z_n)}^2(z), ~~~z\in \D.  $$
From the fact
$$f_{1,z_n}(\varphi(z_n))=f_{1,z_n}(\psi(z_n))=f_{1,z_n}'(\psi(z_n))=0,$$
we get
\begr
&&\|(uC_{\varphi}-vC_{\psi})f_{1,z_n}\|_{\B}\nonumber\\
&\geq& (1-|z_n|^2)|u'(z_n)f_{1,z_n}(\varphi(z_n))-v'(z_n)f_{1,z_n}(\psi(z_n))\nonumber\\
&~&+u(z_n)f_{1,z_n}'(\varphi(z_n))\varphi'(z_n)-v(z_n)f_{1,z_n}'(\psi(w))\psi'(z_n)|\nonumber\\
&=&|u(z_n)\varphi^{\#}(z_n)\varphi(z_n)|\rho(z_n)^2.\nonumber
\endr
Thus,
\begr
\lim_{|\varphi(z_n)| \to 1}|u(z_n)\varphi^{\#}(z_n)|\rho(z_n)^2
&\leq& \lim_{|\varphi(z_n)| \to 1}\|(uC_{\varphi}-vC_{\psi})f_{1,z_n}\|_{\B}\nonumber\\
&\leq& \lim_{|a| \to 1}\sup_{b\in\D}\|(uC_{\varphi}-vC_{\psi})\left((\sigma^2_{a}-a\sigma_{a})\sigma_{b}^2\right)\|_{\B}.\nonumber\\
\nonumber
\endr
So,
$$
\lim_{n \to \infty}|u(z_n)\varphi^{\#}(z_n)|\rho(z_n)^2 = 0. \eqno{(2.1)}$$
Since $|\varphi(z_n)| \to 1$ and $|u(z_n)\varphi^{\#}(z_n)| \nrightarrow 0$, we obtain $\rho(z_n) \to 0.$
Put  $$g_{1,z_n}(z) = (\sigma_{\varphi(z_n)}(z)-\varphi(z_n))\sigma_{\psi(z_n)}^2(z), ~~~z\in \D. $$
Then we have
\begr
&&\|(uC_{\varphi}-vC_{\psi})g_{1,z_n}\|_{\B}\nonumber\\
&\geq& \bigg|(1-|z_n|^2)u'(z_n)\varphi(z_n)\tau(z_n)^2\nonumber\\
&~&-u(z_n)\varphi^{\#}(z_n)\bigg(\tau(z_n)^2-2\varphi(z_n)\tau(z_n)\frac{(1-|\varphi(z_n)|^2)(1-|\psi(z_n)|^2)}{(1-\overline{\psi(z_n)}\varphi(z_n))^2}\bigg)\bigg|\nonumber\\
&\geq& 2|u(z_n)\varphi^{\#}(z_n)\varphi(z_n)|\rho(z_n)(1-\rho(z_n)^2)\nonumber\\
&~&-(1-|z_n|^2)|u'(z_n)\varphi(z_n)|\rho(z_n)^2-|u(z_n)\varphi^{\#}(z_n)|\rho(z_n)^2\nonumber\\
&\geq& |u(z_n)\varphi^{\#}(z_n)|(2|\varphi(z_n)|\rho(z_n)-\rho(z_n)^2-2|\varphi(z_n)|\rho(z_n)^3)\nonumber\\
&~&-(1-|z_n|^2)|u'(z_n)\varphi(z_n)|\rho(z_n)^2,\nonumber
\endr
where
$$\tau(z_n) = \frac{\psi(z_n)-\varphi(z_n)}{1-\overline{\psi(z_n)}\varphi(z_n)}.$$
Since
\begr
\lim_{|a|\to1}\sup_{b\in\D}\|(uC_{\varphi}-vC_{\psi})\left((\sigma_{a}-a)\sigma_b\right)\|_{\B}
\geq\|(uC_{\varphi}-vC_{\psi})g_{1,z_n}\|_{\B},\nonumber
\endr
using condition (B) and (2.1), and noting that $|\varphi(z_n)|\to1$, we get
$$|u(z_n)\varphi^{\#}(z_n)|\rho(z_n) \to 0.$$
We get the condition $(vi)$ of Theorem B.

Since $\rho(z_n)\to 0$ for $\{z_n\}\in\Gamma^{\#}_{u,\varphi},$ we have $\psi(z_n)\to1.$
Now,  set
$$f_{2,z_n}(z)= (\sigma^2_{\psi(z_n)}(z)-\psi(z_n)\sigma_{\psi(z_n)}(z))\sigma_{\varphi(z_n)}^2(z), ~~~z\in \D,$$
and
$$g_{2,z_n}(z) = (\sigma_{\psi(z_n)}(z)-\psi(z_n))\sigma_{\varphi(z_n)}^2(z), ~~~z\in \D. $$
We obtain that
\begr
&&\lim_{|\psi(z_n)| \to 1}|v(z_n)\psi^{\#}(z_n)|\rho(z_n)\nonumber\\
&\lesssim& \lim_{|\psi(z_n)| \to 1}\|(uC_{\varphi}-vC_{\psi})f_{2,z_n}\|_{\B}
+\lim_{|\psi(z_n)| \to 1}\|(uC_{\varphi}-vC_{\psi})g_{2,z_n}\|_{\B}\nonumber\\
&\lesssim& \lim_{|a| \to 1}\sup_{b\in\D}\|(uC_{\varphi}-vC_{\psi})\left((\sigma^2_{a}-a\sigma_{a})\sigma_{b}^2\right)\|_{\B}\nonumber\\
&~&+\lim_{|a| \to 1}\sup_{b\in\D}\|(uC_{\varphi}-vC_{\psi})\left((\sigma_{a}-a)\sigma_{b}^2\right)\|_{\B}\nonumber\\
&=&0.\nonumber
\endr
Next, we consider
$$h_{1,z_n}(z) = (\sigma_{\varphi(z_n)}(z)-\varphi(z_n))\sigma_{\psi(z_n)}(z), ~~~z\in \D. $$
Then we have
\begr
&&\|(uC_{\varphi}-vC_{\psi})h_{1,z_n}\|_{\B}\nonumber\\
&\geq& \bigg|-(1-|z_n|^2)u'(z_n)\varphi(z_n)\tau(z_n)-u(z_n)\varphi^{\#}(z_n)\tau(z_n)\nonumber\\
&~&+u(z_n)\varphi^{\#}(z_n)\varphi(z_n)\frac{(1-|\varphi(z_n)|^2)(1-|\psi(z_n)|^2)}{(1-\overline{\psi(z_n)}\varphi(z_n))^2}\nonumber\\
&~&+v(z_n)\psi^{\#}(z_n)\bigg(\frac{\varphi(z_n)-\psi(z_n)}{1-\overline{\varphi(z_n)}\psi(z_n)}-\varphi(z_n)\bigg)\bigg|\nonumber\\
&=& \bigg|-(1-|z_n|^2)u'(z_n)\varphi(z_n)\tau(z_n)-u(z_n)\varphi^{\#}(z_n)\tau(z_n)\nonumber\\
&~&+u(z_n)\varphi^{\#}(z_n)\varphi(z_n)\bigg(\frac{(1-|\varphi(z_n)|^2)(1-|\psi(z_n)|^2)}{(1-\overline{\psi(z_n)}\varphi(z_n))^2}-1\bigg)\nonumber\\
&~&+\bigg(u(z_n)\varphi^{\#}(z_n)-v(z_n)\psi^{\#}(z_n)\bigg)\varphi(z_n)\nonumber\\
&~&+v(z_n)\psi^{\#}(z_n)\bigg(\frac{\varphi(z_n)-\psi(z_n)}{1-\overline{\varphi(z_n)}\psi(z_n)}\bigg)\bigg|.\nonumber
\endr
By the conditions (B) and Lemma 2.1, we get
\begr
&&\lim_{n \to \infty}|u(z_n)\varphi^{\#}(z_n)-v(z_n)\psi^{\#}(z_n)|\nonumber\\
&\lesssim&\lim_{|\varphi(z_n)|\to1}\|(uC_{\varphi}-vC_{\psi})f_{1,z_n}\|_{\B}
+\lim_{|\varphi(z_n)|\to1}\|(uC_{\varphi}-vC_{\psi})g_{1,z_n}\|_{\B}\nonumber\\
&~&+\lim_{|\psi(z_n)|\to1}\|(uC_{\varphi}-vC_{\psi})f_{2,z_n}\|_{\B}
+\lim_{|\psi(z_n)|\to1}\|(uC_{\varphi}-vC_{\psi})g_{2,z_n}\|_{\B}\nonumber\\
&~&+\lim_{|\varphi(z_n)|\to1}\|(uC_{\varphi}-vC_{\psi})h_{z_n}\|_{\B}\nonumber\\
&\lesssim& \lim_{|a| \to 1}\sup_{b\in\D}\|(uC_{\varphi}-vC_{\psi})\left((\sigma^2_{a}-a\sigma_{a})\sigma_{b}^2\right)\|_{\B}\nonumber\\
&~&+\lim_{|a| \to 1}\sup_{b\in\D}\|(uC_{\varphi}-vC_{\psi})\left((\sigma_{a}-a)\sigma_{b}^2\right)\|_{\B}\nonumber\\
&~&+\lim_{|a| \to 1}\sup_{b\in\D}\|(uC_{\varphi}-vC_{\psi})\left((\sigma_{a}-a)\sigma_{b}\right)\|_{\B}.\nonumber
\endr
Thus, we have
$$\lim_{n \to \infty}|u(z_n)\varphi^{\#}(z_n)-v(z_n)\psi^{\#}(z_n)|=0. \eqno{(2.2)}
$$
We get condition $(v)$ of Theorem B.\msk

For $\{z_n\}\in\Gamma^{\#}_{u,\varphi}$, (2.2) implies
$|v(z_n)\psi^{\#}(z_n)| \nrightarrow 0$.
Thus $\Gamma^{\#}_{u,\varphi}\subset\Gamma^{\#}_{v,\psi}.$ The converse inclusion can be shown by the same way.
So, we get $\Gamma^{\#}_{u,\varphi}=\Gamma^{\#}_{v,\psi},$ for the case $\Gamma^{\#}_{u,\varphi}\neq\varnothing$.
By the fact above, it is easy to see that  $\Gamma^{\#}_{u,\varphi}\neq\varnothing$ if and only if  $\Gamma^{\#}_{v,\psi}\neq\varnothing$.
Hence, we get $(iv)$ of Theorem B.\msk

Assume that $\Lambda_{u,\varphi} \neq \varnothing$
and take a sequence $\{z_n\}\subset\Lambda_{u,\varphi}$. By $(iv)$ of Theorem B, we have
$$\lim_{|\varphi(z_n)|\to1}|u(z_n)\varphi^{\#}(z_n)|\rho(z_n)=0.$$
Using $g_{1,z_n}$ again, we have
\begr
&&\|(uC_{\varphi}-vC_{\psi})g_{1,z_n}\|_{\B}\nonumber\\
&\geq& \bigg|(1-|z_n|^2)u'(z_n)\varphi(z_n)\tau(z_n)^2\nonumber\\
&~&-u(z_n)\varphi^{\#}(z_n)\bigg(\tau(z_n)^2-2\varphi(z_n)\tau(z_n)\frac{(1-|\varphi(z_n)|^2)(1-|\psi(z_n)|^2)}{(1-\overline{\psi(z_n)}\varphi(z_n))^2}\bigg)\bigg|\nonumber\\
&\geq&(1-|z_n|^2)|u'(z_n)\varphi(z_n)|\rho(z_n)^2\nonumber\\
&~&- |u(z_n)\varphi^{\#}(z_n)|(2|\varphi(z_n)|\rho(z_n)+\rho(z_n)^2-2|\varphi(z_n)|\rho(z_n)^3).\nonumber
\endr
Therefore
$$\lim_{n\to\infty}(1-|z_n|^2)|u'(z_n)|\rho(z_n)^2=0.$$
For $\{z_n\}\subset\Lambda_{u,\varphi}$, we have
$$\lim_{n\to\infty}|(1-|z_n|^2)|u'(z_n)|\neq0.$$
Hence, we get
$\rho(z_n)\to0$ and
$$\lim_{n\to\infty}(1-|z_n|^2)|u'(z_n)|\rho(z_n)=0.$$
Then we get condition $(iii)$ of Theorem B.\msk

By $\rho(z_n)\to0$, we have $|\psi(z_n)|\to1$. Taking
$$F_{z_n}(z)=\sigma_{\varphi(z_n)}(z)\sigma_{\psi(z_n)}(z)-\varphi(z_n)\psi(z_n),$$
we have
\begr
\|(uC_{\varphi}-vC_{\psi})F_{z_n}\|_{\B}&\geq& (1-|z_n|^2)|u'(z_n)-v'(z_n)||\varphi(z_n)\psi(z_n)|\nonumber\\
&~~&-|u(z_n)\varphi^{\#}(z_n)|\rho(z_n)-|v(z_n)\psi^{\#}(z_n)|\rho(z_n).\nonumber
\endr
Noting that
$$F_{z_n}(z)=h_{1,z_n}(z)+\varphi(z_n)(\sigma_{\psi(z_n)}(z)-\psi(z_n)),$$
we have
\begr
&&\lim_{|\varphi(z_n)|\to1}\|(uC_{\varphi}-vC_{\psi})F_{z_n}\|_{\B}\nonumber\\
&\leq&\lim_{|a|\to1}\sup_{b\in\D}\|(uC_{\varphi}-vC_{\psi})\left((\sigma_a-a)\sigma_b\right)\|_{\B}
+\lim_{|a|\to1}\|(uC_{\varphi}-vC_{\psi})(\sigma_a-a)\|_{\B}.\nonumber
\endr
Therefore, we obtain condition $(ii)$ of Theorem B and $\Lambda_{u,\varphi}=\Lambda_{v,\psi}$ for the case
$\Lambda_{u,\varphi}\neq\varnothing.$ By the fact above, one can check that $\Lambda_{u,\varphi}=\varnothing$
if and only if $\Lambda_{v,\psi}=\varnothing.$ We get condition $(i)$ of Theorem B.
The proof is complete.\\

{\bf Acknowledgement.}  This project was partially supported by NSF of China $($No.11471143 and No.11720101003$)$.

\end{document}